\documentclass[11pt,a4paper]{article}
\usepackage{ifthen,latexsym,amssymb,amsmath,bbm,fixmath}
\usepackage[shortlabels]{enumitem}
\usepackage[nobysame,initials]{amsrefs}

\newcommand{\C}[1]{{\protect\mathcal{#1}}}
\newcommand{\I}[1]{{\mathbbm #1}}
\newcommand{\e}{\varepsilon}

\renewcommand{\mid}{:}
\renewcommand{\ge}{\geqslant}
\renewcommand{\le}{\leqslant}

\newcommand{\beq}[1]{\begin{equation}\label{#1}}
\newcommand{\eeq}{\end{equation}}

\newtheorem{theorem}{Theorem}
\newtheorem{lemma}[theorem]{Lemma}

\newcommand{\bpf}[1][Proof.]{\smallskip\noindent{\it #1} }
\newcommand{\qed}{\nolinebreak\mbox{\hspace{5 true pt}%
  \rule[-0.85 true pt]{3.9 true pt}{8.1 true pt}}}
\newcommand{\epf}{\qed \medskip}

\newtheorem{claim}{Claim}
\newcommand{\cqed}{\nolinebreak\mbox{\hspace{5 true pt}%
  \rule[-0.85 true pt]{2.0 true pt}{8.1 true pt}}}
\newcommand{\bcpf}{\bpf[Proof of Claim.]}
\newcommand{\ecpf}{\cqed \medskip}

\newcommand{\Path}[4]{\ifthenelse{\equal{#1}{}}{P_{#2}(#3,#4)}{P_{#2}(#3,#4,#1)}}
\newcommand{\A}[2]{\overline{#1}(#2)}
\newcommand{\tA}[2]{\widetilde{#1}(#2)}

\newcommand{\Partial}[2]{#1\rightharpoonup #2}
\newcommand{\dom}{\mathrm{dom}}
\newcommand{\Pot}[1]{\Phi(#1)}

\title{Local version of Vizing's theorem for multi-graphs}

\author{Clinton T.\ Conley\\
Department of Mathematical Sciences\\
Carnegie Mellon University\\
Pittsburgh, PA 15213, USA\\
Email: \texttt{clintonc@andrew.cmu.edu}
 \and Jan Greb\'\i k\\
UCLA Mathematics,
 Los Angeles, CA 90095, USA, and\\
 Faculty of Informatics, Masaryk
University, 
602 00 Brno, Czech Republic\\
Email: \texttt{grebikj@math.ucla.edu} \and 
Oleg Pikhurko\\
Mathematics Institute and DIMAP\\
University of Warwick\\
Coventry CV4 7AL, UK\\
Email: \texttt{o.pikhurko@warwick.ac.uk}
}

\begin{document}

\maketitle

\begin{abstract}
 Extending a result of Christiansen, we prove that every mutli-graph $G=(V,E)$ admits a proper edge colouring $\phi:E\to \{1,2,\dots\}$ which is \emph{local}, that is, 
$\phi(e)\le \max\{d(x)+\pi(x),d(y)+\pi(y)\}$ for every edge $e$ with end-points $x,y\in V$, where $d(z)$ (resp.\ $\pi(z)$) denotes the degree of a vertex $z$ (resp.\ the maximum edge multiplicity at $z$). This is derived from a local version of the Fan Equation. \end{abstract}

\section{Introduction}

Edge-colouring is an important and active area of graph theory; for an overview see, for example, the book by Stiebitz, Scheide, Toft and Favrholdt~\cite{StiebitzScheideToftFavrholdt:gec}. 
One of the key results here is the remarkable theorem of Vizing~\cite{Vizing64},  proved independently by Gupta~\cite{Gupta66}, that every multi-graph $G$ admits a proper edge colouring with at most $\Delta(G)+\pi(G)$ colours, where $\Delta(G)$ and $\pi(G)$ denote respectively the maximum degree and the maximum edge multiplicity of~$G$, and an edge colouring is called \emph{proper} if no two different edges sharing a vertex get the same colour.

Some local versions of edge colouring problems (where a possible colour of an edge $xy$ depends on local information such as the degrees of $x$ and $y$ rather than the global maximum degree) were introduced already in the influential paper of Erd\H{o}s, Rubin and Taylor~\cite{ErdosRubinTaylor80}. 
One of the strongest asymptotic results here is by Bonamy, Delcourt, Lang and Postle \cite{BonamyDelcourtLangPostle24} 
giving the following local list version of Vizing's theorem: for every $\e>0$ if the maximum degree $\Delta(G)$ of a graph $G$ is sufficiently large, the minimum degree $\delta(G)$ is at least  $\ln^{25} \Delta(G)$ and we have an assignment of colour lists to edges $e\mapsto L(e)$ so that each edge $\{x,y\}$ gets at least $(1+\e)\max\{d(x),d(y)\}$ colours then there is a proper edge colouring $\phi$ of $G$ with $\phi(e)\in L(e)$ for each edge $e$, where $d(x)$ denotes the degree of a vertex $x$.

Resolving a conjecture posed in \cite{BonamyDelcourtLangPostle24}, Christiansen~\cite{Christiansen22,Christiansen23} proved that every simple graph $G$ admits a proper edge colouring into $\I N:=\{1,2,3,\dots\}$ which is \emph{local}, meaning that the colour of any edge $\{x,y\}$ is at most $\max\{d(x),d(y)\}+1$. In the special case when the graph $G$ is biparite, the stronger conclusion that the colour of $\{x,y\}$ is at most $\max\{d(x),d(y)\}$ (which may be called the \emph{local K\"onig Theorem}) follows from the more general results of Borodin, Kostochka and Woodall~\cite{BorodinKostochkaWoodall97}.

The purpose of this note is to observe that the proof of Christiansen extends to multi-graphs. The following definition seems to be the ``right'' one in this context: let us call an $\I N$-valued edge colouring $\phi$ of a multi-graph \emph{local} if for every edge $e$ with endpoints $x$ and $y$ it holds that
 \beq{eq:main}
 \phi(e)\le \max\{ d(x)+\pi(x),\, d(y)+\pi(y)\},
\eeq
 where $\pi(z)$ denotes the maximum edge multiplicity at a vertex~$z$. Here is a local version of Vizing's theorem for multi-graphs.

\begin{theorem}\label{th:main} Every multi-graph admits a proper local edge colouring.
\end{theorem}

The standard example (a 3-vertex multi-graph where every pair has the same multiplicity) shows that one cannot decrease the right-hand size of~\eqref{eq:main} in general. Note that, even without the extra locality restriction, we do not know the minimum number of colours as a function of $\Delta$ and $\pi$ that suffices for proper edge colouring of every multi-graph with maximum degree $\Delta$ and maximum edge multiplicity $\pi$; see Scheide and Stiebitz~\cite{ScheideStiebitz12} for our current knowledge on this question.

The original Vizing's theorem as well as some other edge colouring results can be derived from the so-called Fan Equation of Vizing~\cite{Vizing65} (which is discussed in detail in~\cite[Section~2]{StiebitzScheideToftFavrholdt:gec}). 
Here we present a local version of the Fan Equation (which, in this more general setting, becomes an inequality, see Theorem~\ref{th:FanEq}) and derive Theorem~\ref{th:main} from it. 
As another consequence  of Theorem~\ref{th:FanEq}, we have the following result.

\begin{theorem}\label{th:Delta} Let $G=(V,E)$ be a simple graph and $k\in \I N$. If the maximum degree of $G$ is at most $k$ and the set of vertices of degree exactly $k$ spans no cycle in $G$, then there is a proper local colouring $E\to \{1,\dots,k\}$.\end{theorem}

\section{Notation and preliminaries}

By $\I N:=\{1,2,\dots\}$ we denote the set of positive integers. For integers $k\ge \ell\ge 0$, we denote $[\ell,k]:=\{\ell,\dots,k\}$ and $[k]:=\{1,\dots,k\}$.

Let $G=(V,E,\i)$ be a \emph{multi-graph}, that is, $V$ and $E$ are the sets of \emph{vertices} and \emph{edges} respectively, and $\i$ is a (not necessarily injective) function from $E$ to ${V\choose 2}$, the set of unordered pairs of vertices. 
The \emph{(edge) multiplicity} $\pi(x,y):=|\i^{-1}(\{x,y\})|$ of a pair $\{x,y\}\in{V\choose 2}$ is the number of edges whose end-points are $x$ and $y$. For $x,y\in V$ and $e\in E$, we write $x\in e$ (resp.\ $e=xy$) to mean that $x\in \i(e)$ (resp.\ $\i(e)=\{x,y\}$). 

The \emph{degree}  of $x\in V$ is $d(x):=\sum_{y\in V\setminus \{x\}} \pi(x,y)$ which is the number of edges incident to $x$. Also, let $\pi(x):=\max_{y\in V\setminus \{x\}} \pi(x,y)$ denote  the maximum multiplicity at~$x$.

As only edges will be coloured, we will usually say just ``colouring'' instead of ``edge colouring''.
A \emph{partial colouring} of $G$ is a function $\phi:\dom(\phi)\to\I N$ where the domain $\dom(\phi)$ of $\phi$ is a subset of~$E$. We may also write this as $\phi:\Partial{E}{\I N}$. If $\dom(\phi)=E$ (that is, every edge is coloured) then we call $\phi$ a \emph{colouring} (and write $\phi:E\to\I N$).
A (partial) colouring $\phi$ is \emph{proper} if no two distinct edges sharing at least one vertex get the same colour. It is called \emph{local} if
\beq{eq:Locality} 
 \phi(e)\le \max\big\{d(x)+\pi(x)\mid x\in e\big\},\quad\mbox{for every $e\in \dom(\phi)$},
 \eeq  
 that is, for every coloured edge $e$ 
 there is $x\in e$ with $\phi(e)\le d(x)+\pi(x)$.

Given a partial colouring $\phi$, a \emph{$\phi$-chain} is a sequence of distinct edges $C=(e_1,\dots,e_p)$ such that $\{e_1,\dots,e_p\}\cap \dom(\phi)=\{e_2,\dots,e_p\}$ (i.e.\ $e_1$ is the only uncoloured edge) and, for every $i\in [p-1]$, the edges $e_i$ and $e_{i+1}$ share exactly one vertex. The \emph{$C$-shift} (or \emph{shift along $C$}) of $\phi$ is the partial colouring $\phi'$ with $\dom(\phi')=(\dom(\phi)\setminus\{e_p\})\cup\{e_1\}$ which coincides with $\phi$, except $\phi'(e_i):=\phi(e_{i+1})$ for $i\in [p-1]$ (while $e_p$ is uncoloured under~$\phi'$). Informally speaking, we shift colours one step down along~$C$.
A $\phi$-chain $(e,e')$ with $e=xy$ and $e'=xz$ (thus, $e$ is uncoloured, $e'$ is coloured and $y\not=z$) is called \emph{$\phi$-safe} if  $\phi(e')\le d(y)+\pi(y)$ and no edge incident to $y$ has colour $\phi(e')$.
A $\phi$-chain $C=(e_1,\dots,e_m)$ with $m\ge 3$ is \emph{$\phi$-safe} if, for every $i\in [m-1]$,
the $\phi_i$-chain $(e_{i},e_{i+1})$ is $\phi_i$-safe, where $\phi_i$ denotes the $(e_1,\dots,e_{i})$-shift of~$\phi$. In other words, $C$ is $\phi$-safe if each individual 1-edge step of the $C$-shift is safe with respect to the current colouring. Note that the sequence $C=(e_1)$ made of a single edge $e_1$ uncoloured under $\phi$ is a $\phi$-chain (and its shift keeps the partial colouring $\phi$ unchanged); also, let us agree that every such single-edge chain is $\phi$-safe.
When the colouring $\phi$ is understood, we may say just ``chain'' instead of ``$\phi$-chain'', etc.

For a vertex $x\in V$ and distinct colours $\alpha,\beta\in\I N$ with $\beta$ not present at $x$, let
$\Path{\phi}{x}{\alpha}{\beta}$ 
denote the sequence of edges on the maximal $(\alpha,\beta)$-bichromatic path that starts at the vertex~$x$. Thus the values of $\phi$ on $\Path{\phi}{x}{\alpha}{\beta}$ alternate between $\alpha$ and $\beta$, starting with $\alpha$;
if the colour $\alpha$ is not present at $x$, then $\Path{\phi}{x}{\alpha}{\beta}$ is the empty sequence. Note that we exclude the case that $\Path{\phi}{x}{\alpha}{\beta}$ is a double edge by requiring that $\beta$ is missing at $x$ in this definition.

In some proofs, we will be given an integer $k$ and  the set of possible colours will be restricted to~$[k]$. Then we will use the following definitions that implicitly depend on $k$. 
Given $\phi:\Partial{E}{[k]}$, we let for a vertex $x\in V$ 
$$
\A{\phi}{x}:=[k]\setminus \{\phi(e)\mid e\ni x\},
$$ 
to be the set of colours \emph{missing}~at $x$, and
$$
 \tA{\phi}{x}:=\A{\phi}{x}\cap [d(x)+\pi(x)],
 $$
 to be the set of colours which are \emph{safe} at~$x$. Thus a colour is safe at $x$ if it is missing at $x$ and is at most $d(x)+\pi(x)$. 
 Also, the \emph{potential (function)} of $\phi:\Partial{E}{[k]}$ is 
$$
\Pot{\phi}:=\sum_{x\in V} \left|\tA{\phi}{x}\right|.
$$

\begin{lemma}\label{lm:GoodShift} Let $G=(V,E,\i)$ be a multigraph, $k\in\I N$ be a positive integer,  $\phi:\Partial{E}{[k]}$ be a proper local partial colouring, and $\phi'$ be the shift of $\phi$ along a $\phi$-safe chain $(e,e')$, with  $e=xy$ and $e'=xz$. Then $\phi'$ is a proper local partial colouring and $\Pot{\phi'}\le \Pot{\phi}$. Moreover, we have equality if and only if $\phi(e')\le d(z)+\pi(z)$.\end{lemma}

\bpf The edge $e=xy$ gets coloured during the shift. The new colouring $\phi'$ is proper because the set of colours at $x$ does not change while the new colour on $e=xy$ was missing
at~$y$. Also,  the locality condition for $e$ holds from the $y$-side, that is, $\phi'(e)=\phi(e')\le d(y)+\pi(y)$.

Consider the difference $\Pot{\phi'}-\Pot{\phi}$. Since the contributions of a vertex that sees the same sets of colours on incident edges in $\phi$ and $\phi'$ cancel each other, we have to look at $y$ and $z$ only. The contribution of $y$ to $\Pot{\phi'}-\Pot{\phi}$ is exactly $-1$, because the shift reduces the number of safe colours at $y$ by 1 (namely, $\phi(e')$ is now gone from this list).
The contribution of $z$ to $\Pot{\phi'}-\Pot{\phi}$ is at most $1$, and it is equal to $1$ if and only if the colour moved out from $z$ is in $[d(z)+\pi(z)]$, that is, if and only if $\phi(e')\le d(z)+\pi(z)$. This finishes the proof of the lemma.\epf


\section{Multi-fans}

Let $k\in\I N$ and a multigraph $G=(V,E,\i)$ be given. Let $\phi:\Partial{E}{[k]}$ be any proper local partial colouring. Let $e\in E$ be an uncoloured edge and let $x\in e$.

A \emph{multi-fan at $x$ with respect to $e$ and $\phi$ (and $k$)} is a sequence $F=(e_1,y_1,\dots,e_p,y_p)$ with $p\ge 1$ consisting of edges $e_1,\dots,e_p$ and vertices $y_1,\dots,y_p$ satisfying the following conditions.
\begin{enumerate}[(F1)]
 \item\label{it:F1} The edges $e_1,\dots,e_p$ are distinct, $e_1=e$, and $e_i=xy_i$ for $i\in [p]$.
 \item\label{it:F2} For every $i\in [2,p]$ there is $j\in [i-1]$ such that $\phi(e_i)\in\tA{\phi}{y_j}$; in particular, $e_i$ is coloured by~$\phi$.
 \end{enumerate}
 Also, we denote $V(F):=\{y_1,\dots,y_p\}$. Note that $V(F)$ does not include~$x$.

Usually, the Fan Equation is stated for an \emph{edge critical} multi-graph (when a desired colouring of the whole edge set exists when we remove any edge). Having in mind possible algorithmic and descriptive set theory applications, we state a version where the multi-graph $G$ we want to colour need not be edge critical and the presented result can be used to ``improve'' the current partial colouring using chains of special kind.
With this in mind, we make the following definitions. 

A chain $C=(e_1,\dots,e_p)$ is \emph{improving} if it is $\phi$-safe and the $C$-shift $\phi'$ of $\phi$ has strictly smaller potential, or there is a way to extend $\phi'$ to $e_p$ keeping $\phi'$ proper and local (in particular, strictly decreasing the potential). 
Note that the new proper local colouring $\phi'$ satisfies $\Pot{\phi'}<\Pot{\phi}$ in either case.
We say that a chain $C$ is \emph{Vizing} if it starts with some edges $f_1,\dots,f_q$, $q\ge 1$, all containing the same vertex $z$ and then continues with a (possibly empty) initial segment of a bichromatic path starting at~$y$, where $y\not=z$ is the other endpoint of~$f_q$. (Some further restrictions can be put on possible chains arising from the proof of Theorem~\ref{th:FanEq} but we would like to keep this definition simple and short.)
Also, for $z\in V$ let $d_\phi(z):=\left|\{e\in\dom(\phi)\mid e\ni z\}\right|$ denote the number of edges at $z$ that are assigned a colour by~$\phi$.

\begin{theorem}[Local Fan Inequality]
\label{th:FanEq}
Let $k\ge 2$, let $G=(V,E,\i)$ be a multi-graph with $\Delta(G)\le k$, and let $\phi:\Partial{E}{[k]}$ be a partial proper local colouring that admits no improving Vizing chain. 
Suppose that $e\in E$ is uncoloured. Let $x\in e$ and let $F=(e_1,y_1,\dots,e_p,y_p)$ be a multi-fan at $x$ with respect to $e$ and $\phi$. Then all of the following properties hold.
 \begin{enumerate}[(a)]
  \item\label{it:0} For every $i\in [p]$ there are $m\ge 0$ and a sequence $i_0>i_1>\dots>i_m$ such that $i_0=i$, $i_m=1$ and $(e_{i_m},e_{i_{m-1}},\dots,e_{i_0})$ is a $\phi$-safe chain.
 \item\label{it:a} For every $i\in [p]$, we have
  \beq{eq:a}
  \A{\phi}{x}\cap \A{\phi}{y_i}\cap \left[\,\max\{d(x)+\pi(x),\,d(y_i)+\pi(y_i)\}\,\right]=\emptyset.
  \eeq
  \item\label{it:b} For all choices of $i\in [p]$, $\alpha\in\tA{\phi}{x}$, and $\beta\in \tA{\phi}{y_i}$, the $(\alpha,\beta)$-bichromatic path  $\Path{\phi}{y_i}{\alpha}{\beta}$ from $y_i$ ends in~$x$.
  \item\label{it:c} For every $i,j\in [p]$, if $y_i\not=y_j$ then
  \beq{eq:c}
  \tA{\phi}{y_i}\cap \tA{\phi}{y_j}=\emptyset.
  \eeq
  \item\label{it:d} If the multi-fan $F$ is maximal then $|\{y_1,\dots,y_p\}|\ge 2$ and 
  \beq{eq:d}
   \sum_{z\in V(F)}\left(\pi_F(x,z)-\min\big\{k,d(z)+\pi(z)\big\}+d_\phi(z)\right)\ge 1.
   \eeq
   (Here $\pi_F$ denotes the edge multiplicity in the multi-graph $F$, while $\pi=\pi_G$ is taken with respect to the whole multi-graph~$G$.)
 \end{enumerate}
\end{theorem}

\bpf To prove Part~\ref{it:0} for any given $i\in [p]$, we construct a required sequence as follows. Initially, let $i_0:=i$ and $m:=0$. If the current $i_m$ is equal to 1 then stop. Otherwise let $i_{m+1}\in [i_m-1]$ be any index satisfying~\ref{it:F2} for $i_m$ (that is, $\phi(e_{i_{m}})\in \tA{\phi}{e_{i_{m+1}}}$), increase $m$ by $1$, and repeat. Let $i_0>i_1>\dots>i_m$ be the final sequence (with $i_0=i$, $i_m=1$, and $m\ge 0$). Clearly, $(e_{i_m},\dots,e_{i_0})$ is a chain.

Although this is intuitively obvious, let us formally check that the chain $(e_{i_m},\dots,e_{i_0})$ is $\phi$-safe. This vacuously holds if $m=0$, so assume that $m\ge 1$. Let $\phi_1:=\phi$ and, inductively for $j\in [m]$, let $\phi_{j+1}$ be the $(e_{i_{m-j+1}},e_{i_{m-j}})$-shift of $\phi_j$; equivalently, $\phi_{j+1}$ is obtained from $\phi$ by shifting along $(e_{i_m},\dots,e_{i_{m-j}})$. We show by induction on $j\in [m]$ that the $\phi_j$-chain $(e_{i_{m-j+1}},e_{i_{m-j}})$ is $\phi_j$-safe  and the partial colouring $\phi_{j+1}$ is proper and local. So take any $j\in [m]$. We know that $\phi_j$ is proper and local (by induction if $j\ge 2$ and by $\phi_1=\phi$ if $j=1$). 
Note that $\phi_j(e_{i_{m-j}})=\phi(e_{i_{m-j}})$ (since, by~\ref{it:F1}, the edge $e_{i_{m-j}}$ is distinct from any edge whose colour changes when we construct $\phi_j$ from $\phi_1=\phi$).  By the choice of $i_{m-j+1}$ (that is, by~\ref{it:F2}), $\phi(e_{i_{m-j}})$ is in $\tA{\phi}{y_{i_{m-j+1}}}$. Since all colours at $x$ are pairwise distinct, the colour  $\phi(e_{i_{m-j}})$ cannot appear at $y_{i_{m-j+1}}$ when we pass from $\phi_1=\phi$ to $\phi_j$. Thus $\phi_j(e_{i_{m-j}})=\phi(e_{i_{m-j}})$ is also in $\tA{\phi_j}{y_{i_{m-j+1}}}$ and 
the $\phi_j$-chain $(e_{i_{m-j+1}},e_{i_{m-j}})$ is $\phi_j$-safe. It follows that the partial colouring $\phi_{j+1}$ is proper and local by Lemma~\ref{lm:GoodShift}, as required.

Suppose that Part~\ref{it:a} is false. Let $i\in [p]$ be the smallest index violating it. Let $\alpha$ be a colour present in the left-hand side of~\eqref{eq:a}. Let $(i_0,\dots,i_m)$ be the sequence returned by Part~\ref{it:0} for this index~$i$.
Let $\phi'$ be the shift of $\phi$ along $(e_{i_m},e_{i_{m-1}},\dots,e_{i_0})$. By Lemma~\ref{lm:GoodShift}, the new colouring $\phi'$ is proper and local (and satisfies $\Pot{\phi'}\le \Pot{\phi}$). Under $\phi'$, the edge $e_{i}=xy_i$ is uncoloured; also, the colour $\alpha$ is missing at both $x$ and $y_i$. Indeed, $\alpha$ was missing at both $x$ and $y_i$ before the shift (i.e.\ under~$\phi$) while the shift, that affects only edges~at $x$, does not move colour $\alpha$ at all.

Obtain $\phi''$ from $\phi'$ by colouring $e_i$ with the colour~$\alpha$. By above, the new colouring $\phi''$ is still proper. By definition, $\alpha$ is at most $\max\{d(x)+\pi(x),\,d(y_i)+\pi(y_i)\}$, so $\phi''$ is also local. Thus $(e_{i_m},\dots,e_{i_0})$ is an improving Vizing chain for $\phi$. This contradiction proves Part~\ref{it:a}.

Let us turn to Part~\ref{it:b}. Suppose that the claim is false for some vertex $y_i$ and colours $\alpha$ and~$\beta$. Choose the smallest possible such index~$i\in [p]$. 
Recall that $\beta$ is missing at $y_i$, $P:=\Path{\phi}{y_i}{\alpha}{\beta}$ is the $(\alpha,\beta)$-bichromatic path starting at~$y_i$, and we assumed on the contrary that the final endpoint of $P$ is a vertex $x'\not=x$. The vertex set $V(P)$ of $P$ does not contain $x$, as otherwise $x$ would be an endpoint (since $\alpha$ is missing at $x$). It follows that none of the edges incident to $x$ (in particular, none of $e_1,\dots,e_p$) can belong to~$P$. Let $i_0=i>\dots>i_m=1$ be the sequence returned by Part~\ref{it:0} for this $i$ and let $\phi'$ be the shift of $\phi$ along the chain $(e_{i_m},\dots,e_{i_0})$. By Lemma~\ref{lm:GoodShift}, the partial colouring $\phi'$ is proper and local, and satisfies $\Pot{\phi'}\le \Pot{\phi}$. Note that the edge $e_i$ is not coloured by~$\phi'$. 

Let the path $P$ traverse edges $(f_1,\dots,f_\ell)$ and vertices $(u_1,u_2,\dots,u_{\ell+1}$) in this order. Thus $u_1=y_i$ and $f_s=u_{s}u_{s+1}$ for every $s\in[\ell]$. By Part~\ref{it:a}, the colour $\alpha$ is present at $y_i$ under $\phi$ and also under $\phi'$ (since the other endpoint of the colour-$\alpha$ edge at $y$ cannot be~$x$). Thus $\ell\ge 1$. Let $f_0:=e_i$ denote the edge coming before the path $P$ in our shifting procedure. Let $P':=(f_0,\dots,f_\ell)$. It is a path (since $V(P)\not \ni x$) and a $\phi'$-chain (since $f_0\not\in \dom(\phi')$). 
Also, let $u_{0}:=x$ (so that $f_0=u_{0}u_1$) and let $C$ be the concatenation of $(e_{i_m},\dots,e_{i_0})$ and~$P$. Clearly, $C$ is a Vizing chain.

Observe that if we were interested in only proper edge-colourings (that is, not requiring that the extra locality property in~\eqref{eq:Locality} holds), then we could have taken the $P'$-shift of $\phi'$ and assign one of $\alpha$ or $\beta$ to the (now uncoloured) edge $f_\ell$, thus finding a colouring with strictly larger domain than $\phi$. Instead, we proceed as follows. 

Starting with $\phi'$, we iteratively apply  the $(f_s,f_{s+1})$-shifts for $s=0,1,\dots,\ell-1$ as long as each of them is safe. Lemma~\ref{lm:GoodShift} implies that the potential $\Phi$ cannot increase at any step. Thus $\Phi$ stays constant by our assumptions on~$\phi$.

Suppose first that we cannot perform the above shift for some $s\le \ell-1$, that is, before doing the whole chain~$P'$. Denote the obtained colouring $\phi''$. Thus $\phi''$ is the $(f_0,\dots,f_s)$-shift of $\phi'$, where $s$ is the maximum index (assumed to be at most $\ell-1$) such that $(f_0,\dots,f_s)$ is a $\phi'$-safe chain. Note that $s\ge 1$ since the first shift (along $(e_i,f_1)$)  moves the colour $\alpha$ to $e_i=xy_i$ and cannot violate the locality condition at $x$ by $\alpha\le d(x)+\pi(x)$. 
Since the $(f_s,f_{s+1})$-shift is not safe\ in the current colouring $\phi''$ but results in a proper colouring (since we shift along a bichromatic path), we have that 
 \beq{eq:Um-1}
 \gamma:=\phi''(f_{s+1})> d(u_{s})+\pi(u_{s}).
 \eeq
Suppose first that $s\ge 2$. Consider the $(f_{s-2},f_{s-1})$-shift (that is, the penultimate one before we obtained $\phi''$). This shift moves the colour from $f_{s-1}=u_{s-1}u_{s}$ to~$f_{s-2}=u_{s-2}u_{s-1}$, with the moved colour being $\gamma$. Thus by~\eqref{eq:Um-1} and the last claim of Lemma~\ref{lm:GoodShift}, the $(f_{s-2},f_{s-1})$-shift strictly decreased the value of the potential, while the chain $C$ truncated at $f_{s-1}$ is safe and thus improving, a contradiction. So, suppose that $s=1$.
Here when we shift the colouring $\phi'$ along the path $P'$, we perform the first safe shift $(f_0,f_1)$ but we cannot proceed with the second shift, namely along $(f_1,f_2)$. However, then we have that $u_s=y_i$ and $\gamma=\beta$ in~\eqref{eq:Um-1}, contradicting the assumption  of Part~\ref{it:b} that $\beta$ lies in $\tA{\phi}{y_i}\subseteq [d(y_i)+\pi(y_i)]$.

Thus, we can assume that the whole chain $P'$ is $\phi'$-safe (and thus $C$ is $\phi$-safe). We shift $\phi'$ all way along it to obtain~$\phi''$; equivalently, $\phi''$ is the $C$-shift of~$\phi$. The  edge $f_\ell$  is uncoloured by $\phi''$. Let $\gamma$ be the element of $\{\alpha,\beta\}$ different from $\phi''(f_{\ell-1})=\phi'(f_\ell)$. Recall that  $f_{\ell}$ is the last edge of the maximal $\alpha/\beta$-bichromatic path $P$ with respect to $\phi$ that starts with~$y_i$. 
Clearly, the colour $\gamma$ is missing at the penultimate vertex $u_{\ell}$ under $\phi''$. 

Let us show that $\gamma$ cannot occur at $u_{\ell+1}$ under $\phi''$. Suppose otherwise. As  $\gamma$ is not present  at $u_{\ell+1}$ under $\phi$ (by the maximality of $P$) but is present at $u_{\ell+1}$ under $\phi'$, it must be the case that $\gamma=\beta$ and $u_{\ell+1}=y_{i_j}$ for some $j\in [m]$ (and the colour $\beta$ appeared at $y_{i_j}$ because it was shifted to it from $y_{i_{j-1}}$ during the $(e_{i_m},\dots,e_{i_0})$-shift of $\phi$). We conclude that $i_j$ is strictly less than $i=i_0$ and the reversal of the path $P$ is the maximal $\alpha/\beta$-bichromatic path starting at $y_{i_j}$ under~$\phi$. This means that $\Path{\phi}{y_{i_j}}{\alpha}{\beta}$ ends at $y_i\not=x$ and the index $i_j$ contradicts the minimality of $i$. This contradiction proves that $\gamma$ is not present at $u_{\ell+1}$ under~$\phi''$, as claimed.

We conclude that,  under~$\phi''$, the colour $\gamma$ is missing at both $u_{\ell}$ and $u_{\ell+1}$. Next, one can show similarly to above that $\gamma\le d(u_{\ell})+\pi(u_{\ell})$: 
otherwise, for $\ell\ge 2$, the shift for $s=\ell-2$ would strictly decrease~$\Phi$ while, for $\ell=1$, this would contradict the choice of $\beta$ by $\gamma=\beta$ and $u_\ell=y_i$.
Thus if we colour $f_{\ell}$ by~$\gamma$ then we obtain a proper local colouring. Hence, $C$ is an improving Vizing chain. This contradiction proves Part~\ref{it:b}.

Let us turn to Part~\ref{it:c}. Suppose on the contrary there are $i,j\in [p]$ such that $y_i\not=y_j$ and 
there is a colour $\beta\in \tA{\phi}{y_i}\cap \tA{\phi}{y_j}$.
 Let $\alpha$ be any element in $\tA{\phi}{x}$. The last set is non-empty since at most $d(x)-1<k$ edges at $x$ are coloured. Part~\ref{it:b} applies to $i$ (resp.\ $j$) and gives that the $\alpha/\beta$-alternating chain starting at $y_i$ (resp.\ $y_j$) ends at~$x$. This means that, in the multi-graph formed by the edges coloured $\alpha$ or $\beta$, the three distinct vertices $y_i$, $y_j$ and $x$ lie in the same connectivity component and each has degree~$1$. This is contradiction, since the $\{\alpha,\beta\}$-bichromatic multi-graph has maximum degree at most $2$.

For Part~\ref{it:d}, suppose that the multi-fan $F$ is maximal. Since at most $d(y_1)-1$ edges are coloured at $y_1$, there is $\beta\in \tA{\phi}{y_1}$. The colour $\beta$ must be present at $x$ for otherwise, by colouring $e_1$ with colour $\beta$, we obtain a larger proper local colouring and thus $(e_1)$ is an improving chain, a contradiction. Let $e'=xy'$ be the edge of colour~$\beta$. By maximality, the fan $F$ contains the edge~$e'$. Thus $V(F)$ contains $y'\not=y$ and has size at least~$2$, as desired. 
It remains to prove~\eqref{eq:d}. For this, define $\Gamma:=\{\phi(e_2),\dots,\phi(e_p)\}$ and 
$\Gamma':=\bigcup_{z\in V(F)} \tA{\phi}{z}$.
(Recall that $V(F)=\{y_1,\dots,y_p\}$.)

\begin{claim}\label{cl:Gammas} The sets $\Gamma$ and $\Gamma'$ are the same.\end{claim}

\bcpf Let us show that $\Gamma\subseteq\Gamma'$. Take any index $i\in [2,p]$. By~\ref{it:F2}, there is $j\in [i-1]$ such that  $\phi(e_i)\in \tA{\phi}{y_j}$. Thus $\phi(e_i)\in \Gamma'$, as desired.

Conversely, pick any $i\in [p]$ and $\beta\in \tA{\phi}{y_i}$. By Part~\ref{it:a}, we have that $\beta\not\in \A{\phi}{x}$. Thus there is $e'=xz$ with $\phi(e')=\beta$. By the maximality of $F$, there is $j\in [2,p]$ with $e'=e_j$. Thus $\beta=\phi(e_j)$ belongs to $\Gamma$. We conclude that $\Gamma\supseteq \Gamma'$, finishing the proof of the claim.\ecpf 

By Claim~\ref{cl:Gammas} and Part~\ref{it:c}, we have that
\begin{eqnarray*}
 p-1 &=& |\Gamma|\ =\ |\Gamma'|\ =\ \sum_{z\in V(F)} \left| \tA{\phi}{z}\right|
 \ \ge\ \sum_{z\in V(F)} \left(\min\{ k,d(z)+\pi(z)\}-d_\phi(z)\right).
  \end{eqnarray*}
  Also, we have that $p=\sum_{z\in V(F)} \pi_F(x,z)$. Putting these two identities together, we obtain~\eqref{eq:d}. This finishes the proof of the lemma.\epf
  
Note that, in Theorem~\ref{th:FanEq}, if $e$ is the only edge of $G$ not coloured by $\phi$ then $d_\phi(y_i)=d(y_i)$ except $d_\phi(y_1)=d(y_1)-1$ and the Local Fan Inequality simplifies to
 $$
    \sum_{z\in V(F)}\left(\pi_F(x,z)-\min\big\{k-d(z),\pi(z)\big\}\right)\ge2.
 $$

Now we are ready to derive the promised local version of Vizing's theorem.

\bpf[Proof of Theorem~\ref{th:main}.] Let $k:=\max\{d(x)+\pi(x)\mid x\in V(G)\}$. Starting with the empty colouring, iteratively apply improving Vizing chains until none exists. (We stop since the potential strictly decreases each time.)
Suppose that some edge $e$ is not coloured by the final colouring~$\phi$. Let $x\in e$ and let $F=(e_1,y_1,\dots,e_p,y_p)$ be a maximal multi-fan at $x$ with respect to $e$ and~$\phi$. Then the Local Fan Inequality~\eqref{eq:d} holds by Theorem~\ref{th:FanEq}.
Since $k\ge d(z)+\pi(z)$ for every $z\in V(G)$, the inequality states that
$$
 \sum_{z\in V(F)} \left(\pi_F(x,z)-\pi(z)+d_\phi(z)-d(z)\right)\ge 1.
 $$
 This is impossible as each summand in the left-hand side is clearly non-positive.\epf
 
It easy to convert the proof of  Theorem~\ref{th:main} into an algorithm that finds a local colouring of input multi-graph $G=(V,E,\i)$ without isolated vertices with running time polynomial in~$|E|$. The algorithm starts with the empty partial colouring $\phi_\emptyset$ and iteratively finds an improving Vizing chain and changes the current partial colouring accordingly. Each step decreases the potential. Thus we make at most $\Pot{\phi_\emptyset}=O(|E|)$ improvements. Given an uncoloured edge, an improving Vizing chain as in the proof of Theorem~\ref{th:FanEq} can be easily found in polynomial in $|E|$ time.

 Note that if $G$ is a simple graph, then the Local Fan Inequality~\eqref{eq:d} states that
 \beq{eq:d'}
 \sum_{z\in V(F)}\left(d_{\phi}(z)- \min\{k-1,d(z)\}\right)\ge 1.
 \eeq

\bpf[Proof of Theorem~\ref{th:Delta}.] Let $\C C$ be the set of pairs $(\phi,e)$ such that $\phi$ is a proper local partial colouring $\Partial{E(G)}{[k]}$ with the smallest possible potential, while $e\in E$ an edge is not coloured by~$\phi$.  Suppose that the theorem is false. Thus $\C C$ is non-empty.

\begin{claim}\label{cl:Delta} For every $(\phi,e)\in\C C$ and $x\in e$ there are at least two different choice of $f\in E$ such that $f=xy$ for some $y\in V$ with $d(y)=k$ and $(\phi',f)\in\C C$ for some~$\phi'$.\end{claim}

\bcpf Let  $F=(e_1,y_1,\dots,e_p,y_p)$ be a maximal multi-fan at $x$ with respect to $e_1=e$ and~$\phi$. 
By Theorem~\ref{th:FanEq}, the Local Fan Inequality~\eqref{eq:d'} holds. Since $d_\phi(y_1)\le d(y_1)-1$, the contribution of $y_1$ to the left-hand size of~\eqref{eq:d'} is at most 0 if $d(y_1)=k$ and at most $-1$ if $d(y_1)\le k-1$. We conclude by $k\ge d(z)$ for every $z\in V(G)$ that there are at least two indices $i\in [p]$ such that $d(y_i)=k$. 

Let us show that, for each such $i$, the edge $f:=e_i$ satisfies the claim. 
Indeed, let $i_0=i>i_1>\dots> i_m=1$ be the sequence of indices returned by Part~\ref{it:0} of Theorem~\ref{th:FanEq}. Let $\phi'$ be the shift of $\phi$ along the $\phi$-safe sequence $(e_{i_m},e_{i_{m-1}},\dots,e_{i_0})$. By Lemma~\ref{lm:GoodShift}, $(\phi',f)\in\C C$, as required.\ecpf

Start with any $(\phi_1,e_1)\in \C C$ with $e_1=x_1x_2$. By applying Claim~\ref{cl:Delta} (twice) and changing the pair $(\phi_1,e_1)$ we can assume first that $d(x_1)=k$ and then that also $d(x_2)=k$ (and thus $k=\Delta(G)$). Starting with $\phi_1$ and $e_1=x_1x_2$, we inductively construct an infinite sequence $\phi_1,\phi_2,\dots$ and  an infinite path visiting edges $e_1,e_2,\dots$ and vertices $x_1,x_2,x_3,\dots$ in the stated order such that $(\phi_i,e_i)\in\C C$ and $d(x_i)=k$ for every $i\ge 1$ as follows. Let $i\ge 1$ and suppose that we already have $\phi_1,\dots,\phi_i$ and $e_1,\dots,e_i$ as above. Claim~\ref{cl:Delta} when applied to $(\phi_i,e_i)$ with $x=x_{i+1}$ gives at least two different potential choices of the next edge $e_{i+1}= x_{i+1}x_{i+2}$ and the next partial colouring $\phi_{i+1}$ (that is, such that $(\phi_{i+1},e_{i+1})\in \C C$ and $d(x_{i+2})=k$). For at least one of these choices we have $x_{i+2}\not\in\{x_1,\dots,x_i\}$ as otherwise this would create a cycle on vertices of degree~$k$. Thus we can always extend the path by a new edge. However, this contradicts the finiteness of our graph $G$, thus proving Theorem~\ref{th:Delta}.\epf

\section*{Acknowledgements}

The authors would like to thank the anonymous referees for their helpful comments.

Clinton Conley was supported by NSF Grants DMS-1855579 and DMS-2154160.  Jan Greb\'\i k was supported by
Leverhulme Research Project Grant RPG-2018-424. Oleg Pikhurko was supported by  ERC Advanced Grant 101020255 and Leverhulme Research Project Grant RPG-2018-424.



\end{document}